# Sharp weak-type inequalities for differentially subordinated martingales

ADAM OSĘKOWSKI

*Department of Mathematics, Informatics and Mechanics, University of Warsaw, Banacha 2, 02-097 Warsaw, Poland. E-mail: ados@mimuw.edu.pl*

Let $M, N$ be real-valued martingales such that $N$ is differentially subordinate to $M$. The paper contains the proofs of the following weak-type inequalities:

(i) If $M \geq 0$ and $0 < p \leq 1$, then

$$\|N\|_{p,\infty} \leq 2\|M\|_p$$

and the constant is the best possible.

(ii) If $M \geq 0$ and $p \geq 2$, then

$$\|N\|_{p,\infty} \leq \frac{p}{2}(p-1)^{-1/p}\|M\|_p$$

and the constant is the best possible.

(iii) If $1 \leq p \leq 2$ and $M$ and $N$ are orthogonal, then

$$\|N\|_{p,\infty} \leq K_p\|M\|_p,$$

where

$$K_p^p = \frac{1}{\Gamma(p+1)} \cdot \left(\frac{\pi}{2}\right)^{p-1} \cdot \frac{1 + 1/3^2 + 1/5^2 + 1/7^2 + \cdots}{1 - 1/3^{p+1} + 1/5^{p+1} - 1/7^{p+1} + \cdots}.$$

The constant is the best possible.

We also provide related estimates for harmonic functions on Euclidean domains.

## 1. Introduction

The purpose of this paper is to study some sharp estimates for continuous-time martingales. However, to introduce the main concepts and to present the motivations, we will start from the discrete time setting. Let $(\Omega, \mathcal{F}, \mathbb{P})$ be a probability space, filtered by a non-decreasing family $(\mathcal{F}_n)$ of sub-$\sigma$-algebras of $\mathcal{F}$. Let $f = (f_n)$ and $g = (g_n)$ be two real-valued sequences adapted to $(\mathcal{F}_n)$. Let $\mathrm{d}f = (\mathrm{d}f_n)$ and $\mathrm{d}g = (\mathrm{d}g_n)$ be the difference sequences of $f$ and $g$, defined by

$$f_n = \sum_{k=0}^{n} \mathrm{d}f_k, \qquad g_n = \sum_{k=0}^{n} \mathrm{d}g_k, \qquad n = 0, 1, 2, \ldots.$$

Following Burkholder (1989), we say that $g$ is differentially subordinate to $f$ if

$$|\mathrm{d}g_n| \leq |\mathrm{d}f_n|, \qquad n = 0, 1, 2, \ldots \tag{1.1}$$

almost surely. For example, this takes place if $g$ is a transform of $f$ by a predictable real sequence $v = (v_n)$, bounded in absolute value by 1; that is, we have $\mathrm{d}g_n = v_n \mathrm{d}f_n$, $\mathbb{P}(|v_n| \leq 1) = 1$ and $v_n$ is measurable with respect to $\mathcal{F}_{(n-1)\vee 0}$, $n = 0, 1, 2, \ldots$.

Throughout the paper we assume that $f$ and $g$ are $(\mathcal{F}_n)$-martingales. The problem of comparing the sizes of $f$ and $g$ under the assumption of differential subordination has been studied in depth in the literature. For $p \in (0, \infty)$, let

$$\|f\|_p = \sup_n \|f_n\|_p = \sup_n (\mathbb{E}|f_n|^p)^{1/p}$$

and

$$\|f\|_{p,\infty} = \sup_n \|f_n\|_{p,\infty} = \sup_{\lambda > 0} \lambda (\mathbb{P}(f^* \geq \lambda))^{1/p}$$

denote the strong and weak $p$-norms of a martingale. Here $f^* = \sup_n |f_n|$. For $1 < p < \infty$, let $p^* = \max\{p, p/(p-1)\}$. Let us start with the result by Burkholder (1984).

**Theorem 1.1.** *Assume $g$ is differentially subordinate to $f$.*

*We have*

$$\|g\|_{p,\infty} \leq \frac{2}{\Gamma(p+1)} \|f\|_p, \qquad 1 \leq p \leq 2 \tag{1.2}$$

*and*

$$\|g\|_p \leq (p^* - 1) \|f\|_p, \qquad 1 < p < \infty. \tag{1.3}$$

*Both constants $2/\Gamma(p+1)$ and $p^* - 1$ are optimal.*

One can check that neither of the estimates above holds for $p < 1$. The weak-type inequality for the remaining set of parameters $p$ was proved by Suh (2005).

**Theorem 1.2.** *Assume $f$ and $g$ are real-valued and $g$ is differentially subordinate to $f$. Then*

$$\|g\|_{p,\infty} \leq \left( \frac{p^{p-1}}{2} \right)^{1/p} \|f\|_p, \qquad 2 \leq p < \infty. \tag{1.4}$$

*The inequality is sharp.*

If one imposes the additional assumption on the sign of the dominating martingale $f$, the optimal constants change for some values of $p$. Here is one of the main results of Burkholder (1999).

**Theorem 1.3.** *Assume $g$ is differentially subordinate to $f$ and $f$ is non-negative. Then*

$$\|g\|_p \le C_p \|f\|_p, \qquad 1 < p < \infty,$$

*where the optimal constant $C_p$ equals*

$$C_p = \begin{cases} (p-1)^{-1}, & \text{if } 1 < p \le 2, \\ \left[ p\left( \dfrac{p-1}{2} \right) \right]^{1/p}, & \text{if } 2 < p < \infty. \end{cases}$$

Hence the optimal constant in the moment inequalities (1.3) decreases if and only if $2 < p < \infty$. Furthermore, a closer inspection of the proof of (1.2) (see Burkholder (1984), example (4.24), page 657), which shows that the best constant in the inequality (1.2) for non-negative martingale $f$ is still $2/\Gamma(p+1)$. There is a natural question of what can be said if $0 < p < 1$ or $p > 2$. The answer is contained in the following theorem.

**Theorem 1.4.** *Assume $f$ is non-negative, $g$ is real-valued and $g$ is differentially subordinate to $f$. Then*

$$\|g\|_{p,\infty} \le 2\|f\|_p, \qquad 0 < p < 1 \tag{1.5}$$

*and*

$$\|g\|_{p,\infty} \le \frac{p}{2}(p-1)^{-1/p}\|f\|_p, \qquad 2 \le p < \infty. \tag{1.6}$$

*The inequalities are sharp. They are already sharp if $g$ is assumed to be a transform of $f$.*

Now let us turn to the continuous-time setting. Suppose $(\Omega, \mathcal{F}, \mathbb{P})$ is a complete probability space, equipped with a filtration $(\mathcal{F}_t)_{t\ge 0}$, such that $\mathcal{F}_0$ contains all the events of probability 0. Let $M = (M_t)$ and $N = (N_t)$ be two real-valued semimartingales, which have right-continuous paths with limits from the left. The continuous-time extension of the differential subordination, which is due to Bañuelos and Wang (1995) (see also Wang (1995)), can be formulated as follows: The semimartingale $N$ is differentially subordinate to $M$ if the process $([M,M]_t - [N,N]_t)$ is non-negative and non-decreasing. Here $([M,M]_t)$ denotes the quadratic variation process of $M$; see Dellacherie and Meyer (1982). This notion is a generalization of (1.1). To see this, note that if one treats discrete-time sequences $f, g$ as continuous-time processes, then

$$[f,f]_t - [g,g]_t = \sum_{k=0}^{\lfloor t \rfloor} (|\mathrm{d}f_k|^2 - |\mathrm{d}g_k|^2)$$

is non-negative and non-decreasing if and only if (1.1) is valid.

As an example, assume $X$ is a real-valued martingale and $K = (K_s)$ and $H = (H_s)$ are predictable processes such that $|H| \leq |K|$ with probability 1. Let $M$, $N$ be the Itô integrals of $K$, $H$ with respect to $X$; that is,

$$M_t = K_0 X_0 + \int_0^t K_s \, \mathrm{d}X_s, \qquad N_t = H_0 X_0 + \int_0^t H_s \, \mathrm{d}X_s.$$

Then, as

$$[M, M]_t - [N, N]_t = (|K_0|^2 - |H_0|^2)|X_0|^2 + \int_0^t |K_s|^2 - |H_s|^2 \, \mathrm{d}[X, X]_s,$$

we have that $N$ is differentially subordinate to $M$.

All the results above have their counterparts in this new setting. For Theorem 1.1, see the paper by Wang (1995), where a lot of information on transferring inequalities from discrete- to the continuous-time settings can be found. Burkholder's method of proving martingale inequalities involves a construction of a special function, satisfying certain convexity-type properties. Once such a function is found, the continuous-time version follows from Itô's lemma and the smoothing or stopping time argument. For other examples and discussion, see the papers by Bañuelos and Wang (1995) and Suh (2005).

Our approach follows the same pattern. To establish Theorem 1.4, we invent a special function and prove the following stronger result.

**Theorem 1.5.** *Assume $M$ is a non-negative martingale and $N$ is differentially subordinate to $M$. Then*

$$\|N\|_{p,\infty} \leq 2\|M\|_p, \qquad 0 < p < 1, \tag{1.7}$$

*and*

$$\|N\|_{p,\infty} \leq \frac{p}{2}(p-1)^{-1/p}\|M\|_p, \qquad 2 < p < \infty, \tag{1.8}$$

*and the inequalities are sharp.*

We prove the case $p < 1$ in Section 2. Then we deal with the second part of the theorem. As the proof is quite complicated, we divide it into a few steps. First, in Section 3, we show that the constant $p/(2(p-1)^{1/p})$ can not be replaced by a smaller one. Section 4 contains the study of a particular auxiliary differential equation, the solution of which will be needed in Section 5 in order to construct the special function. We complete the proof of Theorem 1.5 in Section 6.

In the second part of the paper we drop the condition $M \geq 0$ and deal with weak-type estimates for differentially subordinated continuous-time martingales under the additional orthogonality assumption. We say that $M$ and $N$ are strongly orthogonal if their covariance process $[M, N]$ is constant with probability 1. In such a case, for convenience, we will skip the word "strongly" and say that $M$ and $N$ are orthogonal.

Our result can be stated as follows:

**Theorem 1.6.** *Assume $M, N$ are real-valued orthogonal martingales with $N$ differentially subordinate to $X$. Then, for $1 \le p \le 2$,*

$$\|N\|_{p,\infty} \le K_p \|M\|_p, \tag{1.9}$$

*where*

$$K_p^p = \frac{1}{\Gamma(p+1)} \cdot \left(\frac{\pi}{2}\right)^{p-1} \cdot \frac{1 + 1/3^2 + 1/5^2 + 1/7^2 + \cdots}{1 - 1/3^{p+1} + 1/5^{p+1} - 1/7^{p+1} + \cdots}. \tag{1.10}$$

*The inequality is sharp.*

The theorem above for the particular case $p = 1$ was proved in Bañuelos and Wang (2000) using the ideas of Choi (1998). Their approach, again based on a construction of a special function, is analytic. In Section 7, we propose a different proof that is more probabilistic in nature.

Finally, the last section of the paper is devoted to related results in harmonic analysis. As exhibited in Burkholder (1991) and Burkholder (1999) (consult also Bañuelos and Wang (1995) for the orthogonal case), the inequalities for differentially subordinated martingales are accompanied by their analogues for harmonic functions on Euclidean domains. Section 8 contains such extensions: the harmonic versions of the inequalities (1.7)–(1.9).

## 2. Theorems 1.4 and 1.5: the case $0 < p < 1$

We start with an auxiliary lemma. Recall that for any semimartingale $X$ there exists a unique continuous local martingale part $X^c$ of $X$ satisfying

$$[X, X]_t = |X_0|^2 + [X^c, X^c]_t + \sum_{0 < s \le t} |\triangle X_s|^2$$

for all $t \ge 0$. Furthermore, $[X^c, X^c] = [X, X]^c$, the pathwise continuous part of $[X, X]$. Here is Lemma 1 from Wang (1995).

**Lemma 2.1.** *The process $Y$ is differentially subordinate to $X$ if and only if $Y^c$ is differentially subordinate to $X^c$, the inequality $|\triangle Y_t| \le |\triangle X_t|$ holds for all $t > 0$ and $|Y_0| \le |X_0|$.*

Now let us introduce the special function $W : \mathbb{R}_+ \times \mathbb{R} \to \mathbb{R}$, constructed in Burkholder (1994) to study the weak-type inequality for non-negative supermartingales. It is given by

$$W(x, y) = \begin{cases} 2x - x^2 + |y|^2, & \text{if } x + |y| \le 1, \\ 1, & \text{if } x + |y| \ge 1. \end{cases}$$

The following functions $\phi:\mathbb{R}_+\times\mathbb{R}\to\mathbb{R}$, $\psi:\mathbb{R}_+\times\mathbb{R}\to\mathbb{R}$ will be needed later:

$$\begin{aligned} &\phi(x,y)=2-2x, \qquad &&\psi(x,y)=2y \qquad &&\text{if } x+|y|\le 1,\\ &\phi(x,y)=0, &&\psi(x,y)=0 &&\text{if } x+|y|>1. \end{aligned}$$

Note that $\phi$ and $\psi$ coincide with the partial derivatives $W_x$, $W_y$ except for the set $\{(x,y):x+|y|=1\}$. It can be shown (see Burkholder (1994), page 1016) that if $x\ge 0$, $x+h\ge 0$, $y,h\in\mathbb{R}$ and $|h|\le|k|$, then

$$W(x+h,y+k)\le W(x,y)+\phi(x,y)h+\psi(x,y)k. \tag{2.1}$$

Furthermore, we have

$$W(x,y)\ge 1_{\{x+|y|\ge 1\}} \tag{2.2}$$

and, if $|y|\le x$, then

$$W(x,y)\le (2x)^p. \tag{2.3}$$

Indeed, the inequality (2.2) is clear; to see (2.3), observe that it suffices to prove it for $|y|=x$ and then the inequality becomes $2x1_{\{2x<1\}}+1_{\{2x\ge 1\}}\le (2x)^p$, which is immediate.

**Proof of the inequality (1.7).** We will prove a stronger statement: for any $\lambda>0$,

$$\lambda^p\mathbb{P}((M+|N|)^*\ge\lambda)\le 2^p\|M\|_p^p.$$

Here, as in the discrete-time case, $X^*=\sup_t|X_t|$. Obviously, we may assume $\lambda=1$. Introduce the stopping time

$$\tau=\inf\{t:M_t+|N_t|\ge 1\}.$$

By Itô's formula,

$$W(M_{\tau\wedge t},N_{\tau\wedge t})=W(M_0,N_0)+I_1+I_2+I_3, \tag{2.4}$$

where

$$\begin{aligned} I_1&=\frac{1}{2}\int_0^{\tau\wedge t}W_{xx}(M_{s-},N_{s-})\,\mathrm{d}[M^c,M^c]_s+2W_{xy}(M_{s-},N_{s-})\,\mathrm{d}[M^c,N^c]_s\\ &\qquad +W_{yy}(M_{s-},N_{s-})\,\mathrm{d}[N^c,N^c]_s=-[M^c,M^c]_{\tau\wedge t}+[N^c,N^c]_{\tau\wedge t},\\ I_2&=\sum_{0<s\le\tau\wedge t}[W(M_s,N_s)-W(M_{s-},N_{s-})-W_x(M_{s-},N_{s-})\Delta M_s-W_y(M_{s-},N_{s-})\Delta N_s],\\ I_3&=\int_0^{\tau\wedge t}W_x(M_{s-},N_{s-})\,\mathrm{d}M_s+\int_0^{\tau\wedge t}W_y(M_{s-},N_{s-})\,\mathrm{d}N_s. \end{aligned}$$

Note that $I_1\le 0$, which is a consequence of Lemma 2.1. Moreover, as

$$W_x(M_{s-},N_{s-})=\phi(W_x(M_{s-},N_{s-})) \quad\text{and}\quad W_y(M_{s-},N_{s-})=\psi(W_x(M_{s-},N_{s-})),$$

we have $I_2 \le 0$: Apply (2.1) to $x = M_{s-}$, $y = N_{s-}$, $h = \Delta M_s$, $k = \Delta N_s$ and observe that $|k| \le |h|$ by Lemma 2.1. Finally, note that $I_3$ is a local martingale. Therefore, there exists a sequence $(T_n)_{n=1}^\infty$ such that $T_n \uparrow \infty$ and, if we replace $t$ with $T_n \wedge t$ $(n = 1, 2, \dots)$ in the expression defining $I_3$, then this expression has zero expectation. Combining this with the previous observations about $I_1$ and $I_2$, we see that (2.4) gives

$$\mathbb{E}W(M_{\tau \wedge T_n \wedge t}, N_{\tau \wedge T_n \wedge t}) \le \mathbb{E}W(M_0, N_0).$$

Now we let $n \to \infty$. As $0 \le W \le 1$, Lebesgue's dominated convergence theorem gives

$$\mathbb{E}W(M_{\tau \wedge t}, N_{\tau \wedge t}) \le \mathbb{E}W(M_0, N_0). \tag{2.5}$$

Apply (2.2) and (2.3) to obtain

$$\mathbb{P}(M_{\tau \wedge t} + |N_{\tau \wedge t}| \ge 1) \le 2^p \mathbb{E}M_0^p = 2^p \|M\|_p^p.$$

To conclude the proof, fix $\varepsilon > 0$ and consider processes $M^\varepsilon = M(1+\varepsilon)$, $N^\varepsilon = N(1+\varepsilon)$. Clearly, $N^\varepsilon$ is differentially subordinate to $M^\varepsilon$, so, if $\eta = \inf\{t : M_t^\varepsilon + |N_t^\varepsilon| \ge 1\}$, we get, by the above argumentation,

$$\begin{aligned}\mathbb{P}((M + |N|)^* \ge 1) = \mathbb{P}((M^\varepsilon + |N^\varepsilon|)^* \ge (1+\varepsilon)) &\le \lim_{t \to \infty} \mathbb{P}(M_{\eta \wedge t}^\varepsilon + |N_{\eta \wedge t}^\varepsilon| \ge 1) \\ &\le 2^p \|M^\varepsilon\|_p^p = 2^p \|M\|_p^p (1+\varepsilon)^p.\end{aligned}$$

As $\varepsilon$ was arbitrary, the proof is complete. □

**Sharpness of (1.5).** Consider the following example: assume the probability space is the interval $[0, 1]$ with Lebesgue measure. Let $f_0 = g_0 \equiv 1/2$ and

$$\mathrm{d}f_1 = -\mathrm{d}g_1 = -\tfrac{1}{2}[0, 3/4] + \tfrac{3}{2}(3/4, 1]$$

and $\mathrm{d}f_n = \mathrm{d}g_n = 0$ for $n \ge 2$. Here and in the next section, we identify a set with its indicator function. Then $|g_1| = 1$ with probability 1 and hence

$$\|g\|_{p,\infty} \ge 1 = 2\|f_0\|_p = 2\|f\|_p,$$

as needed. □

## 3. Sharpness of (1.6) and (1.8)

The optimality of the constant will be proved by constructing an appropriate example. Let $p > 2$ be a fixed number. Let $\delta > 0$, $x \in (0, 1/p)$ be numbers satisfying

$$x\left(1 + \frac{2\delta}{p}\right)^N = \frac{1}{p} \tag{3.1}$$

for the integer $N=N(\delta,x)$. It is clear that we may choose $\delta$ and $x$ to be arbitrarily small.

Consider a two-dimensional Markov martingale $(X_n,Y_n)=(X_n^{x,\delta},Y_n^{x,\delta})$, which is uniquely determined by the following properties:

(i) $X_0=x$, $Y_0=(p-1)x$.

(ii) We have $\mathrm{d}X_n=(-1)^{n+1}\,\mathrm{d}Y_n$ for $n=1,2,\ldots$.

(iii) If $(X_n,Y_n)$ lies on the line $y=(p-1)x$ and $X_n<1/p$, then in the next step it moves either to the line $x=0$ or to the line $y=(p-1+\delta)x/(1+\delta)$, $n=0,1,2,\ldots$.

(iv) If $(X_n,Y_n)$ lies on the line $y=(p-1+\delta)x/(1+\delta)$ then in the next step it moves either to the line $y=(p-1)x$ or to the line $y=(p-2)x/2$, $n=0,1,2,\ldots$.

(v) If $(X_n,Y_n)=(1/p,1-1/p)$ (which happens only if $n=2N$), then $(X_{n+1},Y_{n+1})$ equals either $(0,1)$ or $(2/p,1-2/p)$.

(vi) The states on the line $x=0$ and $y=(p-2)x/2$ are absorbing.

The examplary trajectories are presented on Figure 1.

To be more precise, let the sequence $(p_n)$, $n=0,1,2,\ldots,2N$ be given by

$$p_{2n}=\left(\frac{p-p\delta+4\delta}{(p+2\delta)(1+\delta)}\right)^n,\qquad p_{2n+1}=p_{2n}\cdot\frac{1}{1+\delta},\qquad n=0,1,2,\ldots,N.$$

Let the probability space be the interval $[0,1]$ with Lebesgue measure. Set

$$\begin{aligned}X_0&=x[0,1],\qquad \mathrm{d}X_{2n+1}=\delta X_{2n}[0,p_{2n+1}]-X_{2n}(p_{2n+1},p_{2n}],\\ \mathrm{d}X_{2n+2}&=-(1-2/p)\delta X_{2n}[0,p_{2n+2}]+X_{2n}\left(1+\frac{4}{p}\delta-\delta\right)(p_{2n+2},p_{2n+1}]\end{aligned}\tag{3.2}$$

for $n=0,1,2,\ldots,N-1$, and

$$\mathrm{d}X_{2N+1}=X_{2N}[0,p_{2N}/2]-X_{2N}(p_{2N}/2,p_{2N}].\tag{3.3}$$

Furthermore,

$$Y_0=(p-1)x[0,1],\qquad \mathrm{d}Y_n=(-1)^{n+1}\,\mathrm{d}X_n,\qquad n=0,1,2,\ldots,2N+1.$$

Note that $Y_{2N+1}=1$ on $[0,p_{2N}/2]$ and $0<Y_{2N+1}<1$ on $(p_{2N}/2,1]$, so we have $\mathbb{P}(Y^*\geq 1)=p_{2N}/2$. Furthermore, by (3.2) and (3.3), $X_{2N+1}$ equals 0 on the union of the sets $(p_{2n+1},p_{2n}]$, $n=0,1,\ldots,N$ and the interval $(p_{2N}/2,p_{2N}]$. Moreover, $X_{2N+1}$ is equal to

$$2\left(1+\frac{2\delta}{p}\right)^{n+1}x\qquad\text{on }(p_{2n+2},p_{2n+1}],\qquad n=0,1,\ldots,N-1.$$

Finally, it equals $2/p$ on $[0,p_{2N}/2]$. Hence we may write

$$\mathbb{E}X_{2N+1}^p=\frac{2^p(p-2)\delta}{(p+2\delta)(1+\delta)}\left(1+\frac{2\delta}{p}\right)^p x^p\sum_{k=0}^{N-1}\left(\frac{p-p\delta+4\delta}{(p+2\delta)(1+\delta)}\left(1+\frac{2\delta}{p}\right)^p\right)^k$$

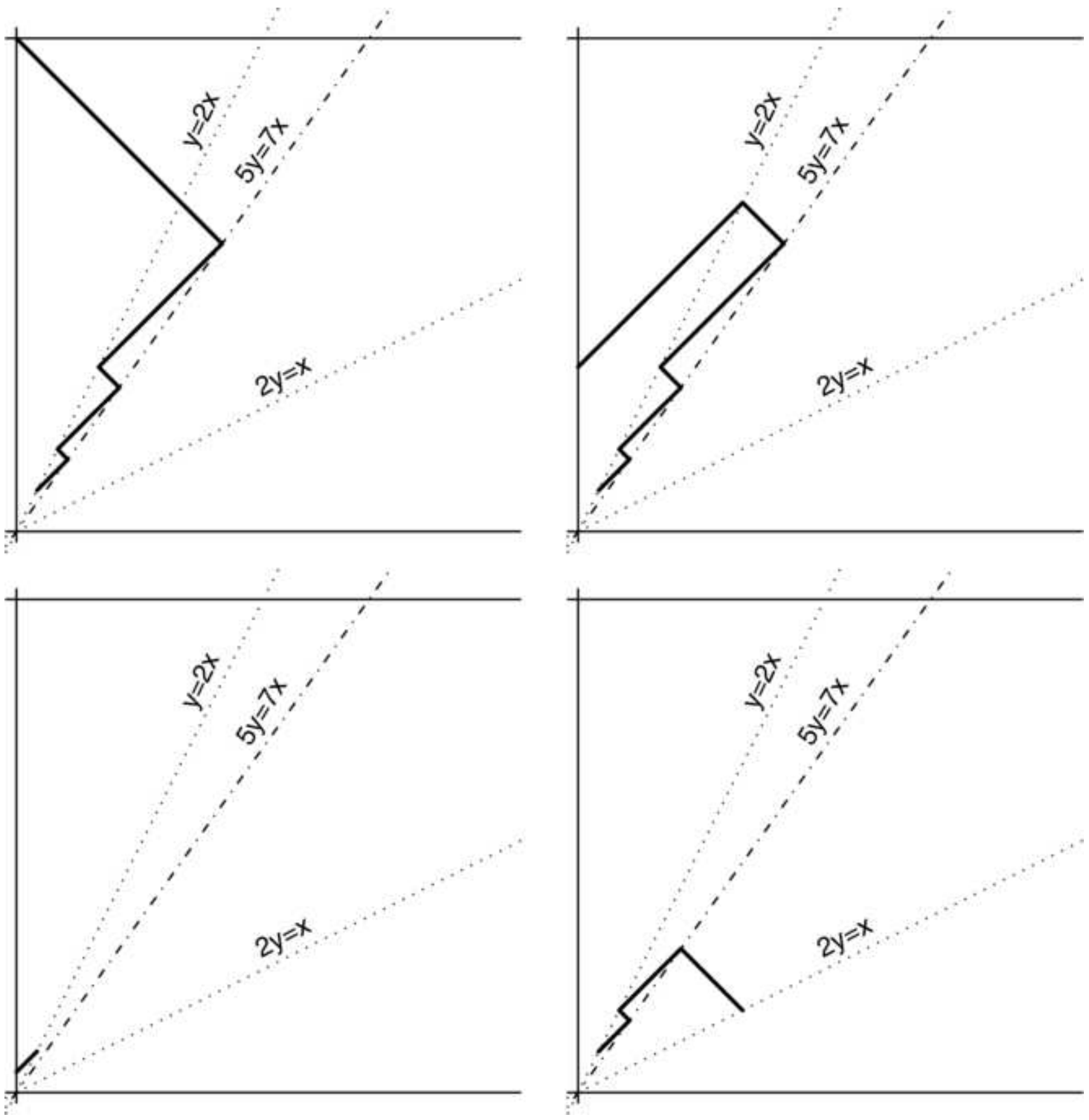

**Figure 1.** Four trajectories of the process $(X, Y)$ corresponding to the parameters $p = 3$, $x = \frac{1}{24}$, $\delta = \frac{3}{2}$ and $N = 3$.

$$
\begin{aligned}
&+ \left(\frac{2}{p}\right)^p \cdot \frac{1}{2}\left(\frac{p - p\delta + 4\delta}{(p+2\delta)(1+\delta)}\right)^N \\
&= \frac{2^p(p-2)\delta}{(p+2\delta)(1+\delta)}\left(1 + \frac{2\delta}{p}\right)^p x^p \\
&\quad \times \frac{((p - p\delta + 4\delta)/((p+2\delta)(1+\delta))(1+2\delta/p)^p)^N - 1}{(p - p\delta + 4\delta)/((p+2\delta)(1+\delta))(1+2\delta/p)^p - 1} \\
&\quad + \left(\frac{2}{p}\right)^p \cdot \frac{1}{2}\left(\frac{p - p\delta + 4\delta}{(p+2\delta)(1+\delta)}\right)^N.
\end{aligned}
$$

Now keep $x$ fixed and let $\delta$ to 0 (so that (3.1) holds, with $N=N(\delta,x)\to\infty$). Then

$$\left(\frac{p-p\delta+4\delta}{(p+2\delta)(1+\delta)}\right)^N=\left(1+\frac{2\delta}{p}\cdot\frac{p(1-p-\delta)}{(p+2\delta)(1+\delta)}\right)^N\to(px)^{p-1}$$

and

$$\frac{\delta}{(p-p\delta+4\delta)/((p+2\delta)(1+\delta))(1+2\delta/p)^p-1}\to\frac{p}{2},$$

so we have

$$\begin{aligned}\mathbb{E}X_{2N+1}^p&\to\frac{2^p(p-2)}{p}x^p\frac{((px)^{p-1}(px)^{-p}-1)p}{2}+\left(\frac{2}{p}\right)^p\cdot\frac{1}{2}(px)^{p-1}\\&=\frac{(2x)^{p-1}}{p}(p-1-p(p-2)x)\end{aligned}$$

and

$$\frac{\mathbb{P}(Y_{2N+1}\geq 1)}{\mathbb{E}X_{2N+1}^p}\to\frac{1/2(px)^{p-1}}{(2x)^{p-1}/p(p-1-p(p-2)x)}=\frac{p^p}{2^p(p-1)-2^pp(p-2)x}. \qquad (3.4)$$

Observe that $Y$ is a transform of $X$ and $|\mathrm{d}X_n|=|\mathrm{d}Y_n|$ for $n\geq 1$. However, $Y$ is not differentially subordinate to $X$ as $Y_0=(p-1)X_0>X_0$. To overcome this difficulty, introduce the processes $Y_n'=(Y_n-(p-2)x)/(1-(p-2)x)$, $X_n'=X_n/(1-(p-2)x)$, $n=0,1,\ldots,2N+1$. Then $Y'$ is a transform of $X'$ by the deterministic sequence $(1,1,-1,1,-1,1,-1,1,\ldots)$ and $\mathbb{P}(Y_{2N+1}'\geq 1)=\mathbb{P}(Y_{2N+1}\geq 1)$. In terms of these new processes, (3.4) reads

$$\frac{\mathbb{P}(Y_{2N+1}'\geq 1)}{\mathbb{E}(X_{2N+1}')^p}\to\frac{p^p(1-x(p-2))^p}{2^p(p-1)-2^pp(p-2)x}$$

and it is clear that the limit can be made arbitrarily close to $p^p/2^p(p-1)$ by choosing $x$ sufficiently small. This proves the sharpness of (1.6) and hence the sharpness of (1.8) as well.

## 4. A differential equation

Let $p>2$ be fixed. The purpose of this section is to study a solution to a certain differential equation. A very similar equation appears in Suh (2005) and our arguments are parallel to those used there. We will show that there exists a function $h:[1,\infty)\to[2/p,\infty)$, which enjoys the following properties:

$$\text{the function } h \text{ is increasing and continuous on } [1,\infty), \qquad (4.1)$$

$$h(t)>t-1 \quad \text{for all } t\geq 1, \qquad (4.2)$$

$$h \text{ is differentiable on } (1,\infty) \quad \text{and} \quad h'(t) = \left(\frac{2}{p}\right)^{p+1} (h(t))^{2-p}(h(t)-t+1)^{-2}, \tag{4.3}$$

$$h(1) = h'(1+) = 2/p. \tag{4.4}$$

The problem above is equivalent to the existence of the function $G:[2/p,\infty) \to [1,\infty)$ satisfying the following properties:

$$\text{the function } G \text{ is increasing and continuous on } [2/p,\infty), \tag{4.5}$$

$$G(t) < t+1 \qquad \text{for all } t \geq 2/p, \tag{4.6}$$

$$G \text{ is differentiable on } (2/p,\infty) \quad \text{and} \quad G'(t) = \left(\frac{p}{2}\right)^{p+1} t^{p-2}(t+1-G(t))^2. \tag{4.7}$$

$$G\left(\frac{2}{p}\right) = 1, \qquad G'\left(\frac{2}{p}+\right) = \frac{p}{2}. \tag{4.8}$$

To see the equivalence, note that if $h$ satisfies (4.1)–(4.4), then $G = h^{-1}$ satisfies (4.5)–(4.8) and if $G$ satisfies (4.5)–(4.8), then $h = G^{-1}$ satisfies (4.1)–(4.4).

As (4.7) has the Riccati form, we can use the transformation

$$k(t) = \exp\left[\int_{2/p}^{t} \left(\frac{p}{2}\right)^{p+1} y^{p-2}(y+1-G(y))\,\mathrm{d}y\right]$$

to obtain the following differential equation for $k$:

$$yk''(y) + (2-p)k'(y) - \left(\frac{p}{2}\right)^{p+1} y^{p-1}k(y) = 0. \tag{4.9}$$

For a fixed $\alpha > -1$, let $I_\alpha$ be the modified Bessel function of the first kind (see Abramowitz and Stegun (1992)). That is,

$$I_\alpha(z) = \sum_{k=0}^{\infty} \frac{(z/2)^{2k+\alpha}}{k!!\Gamma(\alpha+k+1)}$$

and we have

$$z^2 I_\alpha''(z) + zI_\alpha'(z) - (z^2+\alpha^2)I_\alpha(z) = 0.$$

One can check that the functions

$$k_1(t) = t^{(p-1)/2} I_{-(p-1)/p}(z_0)$$

and

$$k_2(t) = t^{(p-1)/2} I_{(p-1)/p}(z_0),$$

where $z_0 = \sqrt{(p/2)^{p-1}t^p}$, are two linearly independent solutions on $(0,\infty)$ to equation (4.9). As the functions $I_\alpha$ are infinitely differentiable on $(0,\infty)$, so are $k_1$ and $k_2$. Let $a_1$, $a_2$ be two numbers such that $k = a_1k_1 + a_2k_2$ satisfies

$$k(2/p) = 1 \quad \text{and} \quad k'(2/p) = \frac{p^2}{4}. \tag{4.10}$$

If one rewrites (4.9) in the form

$$yk''(y) = (p-2)k'(y) + \left(\frac{p}{2}\right)^{p+1} y^{p-1}k(y),$$

then it follows from (4.10) that $k$, $k'$ and $k''$ are strictly positive on $[2/p,\infty)$. Now it is straightforward to check that the function

$$G(t) = t + 1 - \left(\frac{2}{p}\right)^{p+1} \frac{k'(t)}{k(t)t^{p-2}} \tag{4.11}$$

has all the properties (4.5)–(4.8).

We conclude this section with the following lemma:

**Lemma 4.1.** *We have*

$$h'(t) \le 1 \qquad \text{for all } t > 1. \tag{4.12}$$

**Proof.** This is equivalent to $G'(t) \ge 1$ for all $t > 2/p$. We have

$$G''(t) = \left(\frac{p}{2}\right)^{p+1} [(p-2)t^{p-3}(t+1-G(t))^2 + 2t^{p-2}(t+1-G(t))(1-G'(t))],$$

implying that if $G'(t) \le 1$, then $G''(t) > 0$. Now suppose $G'(t_0) < 1$ for some $t_0 > 2/p$. Then, as $G'(2/p+) = p/2 > 1$, there exists $t_1 \in (2/p, t_0)$ such that $G'(t_1) = 1$ and $G'(t) < 1$ for $t \in (t_1, t_0)$. Now by mean value property, for some $t_2 \in (t_1, t_0)$,

$$G''(t_2) = \frac{G'(t_0) - G'(t_1)}{t_0 - t_1} < 0,$$

a contradiction. □

## 5. The special function and its properties

Now we are ready to define the special function $U$. Due to the lack of symmetry with respect to the $y$ axis, this function is much more complicated than the one constructed

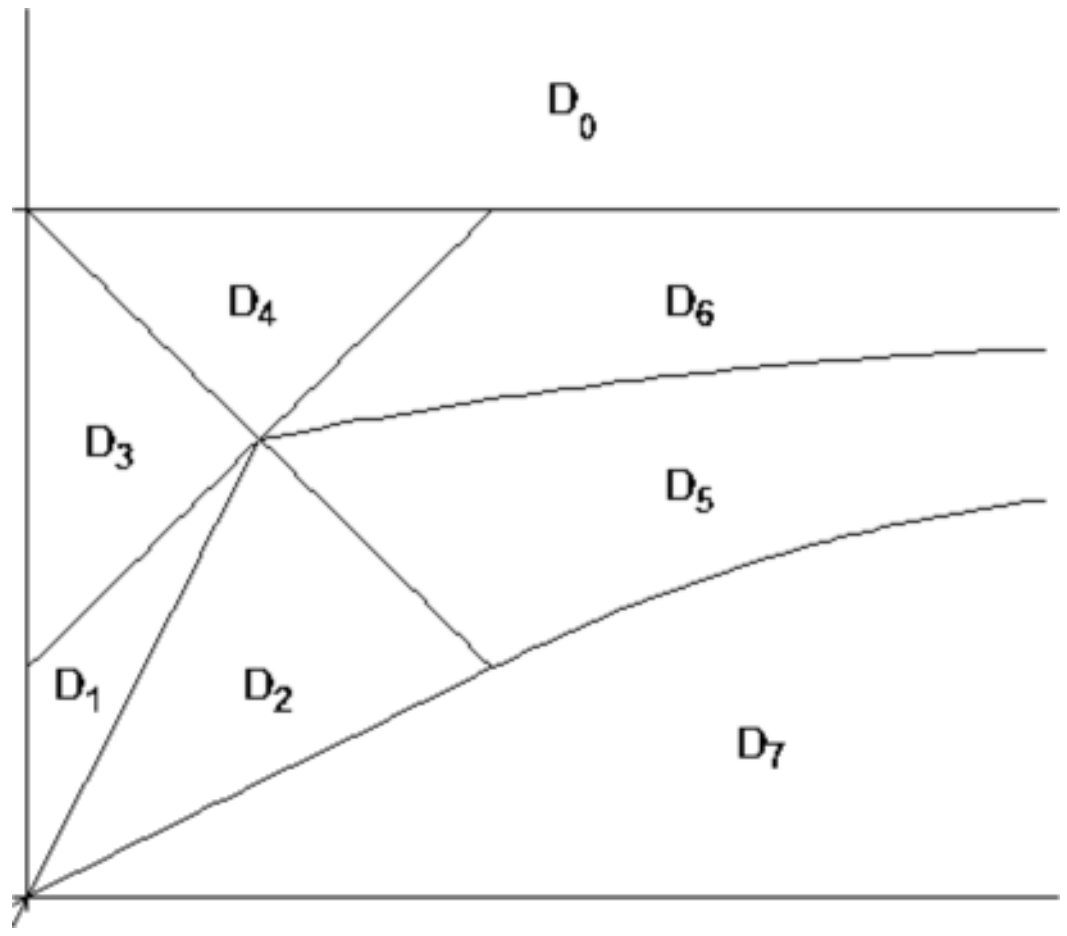

**Figure 2.** The regions $D_0$–$D_7$ for $p=3$, intersected with $\mathbb{R}_+ \times \mathbb{R}_+$.

in Suh (2005). Consider the following subsets of $\mathbb{R}_+ \times \mathbb{R}$.

$$
\begin{aligned}
D_0 &= \{(x,y): |y| \geq 1\},\\
D_1 &= \left\{(x,y): (p-1)x \leq |y| < x+1-\frac{2}{p}\right\},\\
D_2 &= \left\{(x,y): \frac{p-2}{2}x \leq |y| < \min\{1-x,(p-1)x\}\right\},\\
D_3 &= \left\{(x,y): x+1-\frac{2}{p} \leq y < 1-x\right\},\\
D_4 &= \left\{(x,y): \max\left\{1-x, \left(x+1-\frac{2}{p}\right)\right\} \leq |y| < 1\right\},\\
D_5 &= \left\{(x,y): h(x+|y|) \geq x > \frac{-1+h(x+|y|)+x+|y|}{2} \text{ and } x+|y| \geq 1\right\},\\
D_6 &= \left\{(x,y): \frac{1-h(x+|y|)+x+|y|}{2} \leq |y| < \min\left\{\left(x+1-\frac{2}{p}\right),1\right\}\right\},\\
D_7 &= (\mathbb{R}_+ \times \mathbb{R}) \setminus (D_0 \cup D_1 \cup D_2 \cup D_3 \cup D_4 \cup D_5 \cup D_6).
\end{aligned}
$$

See Figure 2 for the case $p=3$.

Let

$$V(x,y)=\begin{cases}1-\dfrac{p^p}{2^p(p-1)}x^p, & \text{on } D_0,\\ -\dfrac{p^p}{2^p(p-1)}x^p, & \text{on } (\mathbb{R}_+\times\mathbb{R})\setminus D_0\end{cases}$$

and let us define $U(x,y)$ by

$$\begin{cases}1-\dfrac{p^p}{2^p(p-1)}x^p, & \text{on } D_0,\\ \dfrac{p^p}{2(p-1)(p-2)^{p-2}}x(|y|-x)^{p-1}, & \text{on } D_1,\\ \dfrac{1}{p-1}(x+|y|)^{p-1}\left[(p-1)|y|-\dfrac{p^2-2p+2}{2}x\right], & \text{on } D_2,\\ \dfrac{x}{2(p-1)(1+x-|y|)}[-(p-2)^2+p^2(|y|-x)], & \text{on } D_3,\\ 1-\dfrac{p^2}{2(p-1)}(1-|y|)-\dfrac{p^p}{2^p(p-1)}(x+|y|-1)(x+1-|y|)^{p-1}, & \text{on } D_4,\\ \dfrac{p^p}{2^p(p-1)}(h(x+|y|))^{p-1}[(p-1)h(x+|y|)-px], & \text{on } D_5,\\ 1-\dfrac{2(1-|y|)}{2+x-|y|-G(x-|y|+1)} & \\ \quad -\dfrac{p^p}{2^p(p-1)}(x-|y|+1)^{p-1}[x-(p-1)(1-|y|)], & \text{on } D_6,\\ -\dfrac{p^p}{2^p(p-1)}x^p, & \text{on } D_7.\end{cases}$$

The properties of $U$ are described in the sequence of the lemmas below.

**Lemma 5.1.** *The function $U$ is continuous on $\mathbb{R}_+\times\mathbb{R}\setminus\{(0,\pm1)\}$ and of class $C^1$ except for the set $\partial D_0\cup(\partial D_3\cap\partial D_4\setminus(\frac{2}{p},\pm\frac{p-2}{p}))$.*

**Proof.** Clearly, $U$ is of class $C^1$ in the interior of each $D_j$, $j=0,1,2,\ldots,7$, so we need to check the properties on the boundaries. By symmetry, we may restrict ourselves to positive $y$'s. Using (4.4), (4.8) and the definitions of the boundaries, the continuity of $U$ can be verified readily. For the second part of the lemma, we calculate the partial derivatives of $U$: we have

$$U_x(x,y)=\begin{cases}-\dfrac{p^{p+1}}{2^p(p-1)}x^{p-1}, & \text{on } D_0^o,\\ \dfrac{p^p}{2(p-1)(p-2)^{p-2}}(y-x)^{p-2}(y-px), & \text{on } D_1^o,\\ \dfrac{p}{2(p-1)}(x+y)^{p-2}[(p-2)y-(p^2-2p+2)x], & \text{on } D_2^o,\\ -\dfrac{p^2}{2(p-1)}+\dfrac{2(1-y)}{(1+x-y)^2}, & \text{on } D_3^o,\\ -\dfrac{p^p}{2^p(p-1)}(x+1-y)^{p-2}[px-(p-2)(1-y)], & \text{on } D_4^o,\\ \dfrac{2(h(x+y)-x)}{(h(x+y)-(x+y)+1)^2}-\dfrac{p^{p+1}}{2^p(p-1)}(h(x+y))^{p-1}, & \text{on } D_5^o,\\ \dfrac{2(1-|y|)}{(2+x-|y|-G(x-|y|+1))^2}-\dfrac{p^{p+1}}{2^p(p-1)}(x-|y|+1)^{p-1}, & \text{on } D_6^o,\\ -\dfrac{p^{p+1}}{2^p(p-1)}x^{p-1}, & \text{on } D_7^o,\end{cases}$$

while

$$U_y(x,y)=\begin{cases}0, & \text{on } D_0^o,\\ \dfrac{p^p}{2(p-2)^{p-2}}x(y-x)^{p-2}, & \text{on } D_1^o,\\ \dfrac{p}{2}(x+y)^{p-2}[2y-(p-2)x], & \text{on } D_2^o,\\ \dfrac{2x}{(1+x-y)^2}, & \text{on } D_3^o,\\ \dfrac{p^2}{2(p-1)}+\dfrac{p^p}{2^p(p-1)}(x+1-y)^{p-2}[(p-2)x-p(1-y)], & \text{on } D_4^o,\\ \dfrac{2(h(x+y)-x)}{(h(x+y)-(x+y)+1)^2}, & \text{on } D_5^o,\\ \dfrac{2}{(2+x-y-G(x-y+1))^2}[1+x-G(x-y+1)], & \text{on } D_6^o,\\ 0, & \text{on } D_7^o.\end{cases}$$

All that is left is to check that the partial derivatives agree on the boundaries. □

**Lemma 5.2.** *For $y\ge 0$, $U_y$ is non-negative.*

**Proof.** This is clear on $D_1^o$, $D_2^o$, $D_3^o$, $D_5^o$ and $D_7^o$. On $D_4$, we have

$$U_{xy}(x,y)=\frac{p^p(p-2)}{2^p}(x+1-y)^{p-3}(x-y)\le 0$$

and, consequently,

$$U_y(x,y) \geq U_y\left(y-1+\frac{2}{p},y\right) = \frac{p^2}{2}\left(y-1+\frac{2}{p}\right) \geq 0.$$

On $D_6$, we have $G(1+x-y) \leq x+y$ (as it is equivalent to $1+x-y \leq h(x+y)$, one of the inequalities defining $D_6$). This can be further bounded from above by $1+x$, which yields the claim. □

The most technical lemma is the following:

**Lemma 5.3.** *Suppose $(x,y)$ belongs to the interior of $D_j$ for some $0 \leq j \leq 7$. Then for any $h$, $k$ we have*

$$U_{xx}(x,y)h^2 + 2U_{xy}(x,y)hk + U_{yy}(x,y)k^2 \leq 0. \tag{5.1}$$

**Proof.** We start with the observation that the inequality holds if $(x,y)$ belongs to $D_0^\circ$ or $D_7^\circ$; indeed, $U_{xy} = U_{yy} = 0$ and $U_{xx} \leq 0$ there. For $(x,y)$ lying in the interior of one of the remaining sets, note that $U$ has the following property: one of the functions $t \mapsto U(x+t,y+t)$, $t \mapsto U(x+t,y-t)$, is linear on some neighbourhood of 0. Hence

$$U_{xx}(x,y) + 2U_{xy}(x,y) + U_{yy}(x,y) = 0 \tag{5.2}$$

or

$$U_{xx}(x,y) - 2U_{xy}(x,y) + U_{yy}(x,y) = 0. \tag{5.3}$$

Now if (5.2) holds, we may write

$$\begin{aligned} &U_{xx}(x,y)h^2 + 2U_{xy}(x,y)hk + U_{yy}(x,y)k^2 \\ &\quad = \frac{U_{xx}(x,y)+U_{yy}(x,y)}{2}(h-k)^2 + \frac{U_{xx}(x,y)-U_{yy}(x,y)}{2}(h^2-k^2), \end{aligned}$$

while if (5.3) is valid, we have

$$\begin{aligned} &U_{xx}(x,y)h^2 + 2U_{xy}(x,y)hk + U_{yy}(x,y)k^2 \\ &\quad = \frac{U_{xx}(x,y)+U_{yy}(x,y)}{2}(h+k)^2 + \frac{U_{xx}(x,y)-U_{yy}(x,y)}{2}(h^2-k^2). \end{aligned}$$

Therefore (5.1) will hold once we have established the inequality

$$U_{xx} \leq -|U_{yy}|. \tag{5.4}$$

As previously, with no loss of generality we may assume $y > 0$.

Straightforward computations show that

$$U_{xx}(x,y) = \begin{cases} \dfrac{p^p}{2(p-2)^{p-2}}(y-x)^{p-3}(px-2y), & \text{on } D_1^o, \\ -p(x+y)^{p-3}\left(\dfrac{p^2-2p+2}{2}x+y\right), & \text{on } D_2^o, \\ -\dfrac{4}{(1+x-y)^3}(1-y), & \text{on } D_3^o, \\ -\dfrac{p^p}{2^p}(x+1-y)^{p-3}(px+(p-4)(y-1)), & \text{on } D_4^o, \\ \dfrac{2}{(h(x+y)-(x+y)+1)^2} & \\ \quad\times\left[-2+h'(x+y)-\dfrac{2(h(x+y)-x)(h'(x+y)-1)}{h(x+y)-(x+y)+1}\right], & \text{on } D_5^o, \\ -\dfrac{4(1-y)(1-G'(x-y+1))}{(2+x-y-G(x-y+1))^3}-\dfrac{p^{p+1}}{2^p}(x-y+1)^{p-2}, & \text{on } D_6^o \end{cases}$$

and

$$U_{yy}(x,y) = \begin{cases} \dfrac{p^p}{2(p-2)^{p-2}}(y-x)^{p-3}\cdot(p-2)x, & \text{on } D_1^o, \\ -p(x+y)^{p-3}\left(\dfrac{p^2-4p+2}{2}x-(p-1)y\right), & \text{on } D_2^o, \\ \dfrac{4}{(1+x-y)^3}x, & \text{on } D_3^o, \\ \dfrac{p^p}{2^p}(x+1-y)^{p-3}(-(p-4)x+p(1-y)), & \text{on } D_4^o, \\ \dfrac{2}{(h(x+y)-(x+y)+1)^2} & \\ \quad\times\left[-h'(x+y)-\dfrac{2(h(x+y)-x)(h'(x+y)-1)}{h(x+y)-(x+y)+1}\right], & \text{on } D_5^o, \\ -\dfrac{4(1-y)(1-G'(x-y+1))}{(2+x-y-G(x-y+1))^3}+\dfrac{2(2-G'(x-y+1))}{(2+x-y-G(x-y+1))^2}, & \text{on } D_6^o. \end{cases}$$

Now let us check (5.4). On $D_1^o$ it is equivalent to

$$px-2y\le -(p-2)x \quad \text{or} \quad y\ge (p-1)x,$$

which follows from the definition of $D_1$. On $D_2^o$, the inequalities $U_{xx}+U_{yy}\le 0$, $U_{xx}-U_{yy}\le 0$ can be transformed to $-p(x+y)\le 0$ and $(p-2)(y-(p-1)x)\le 0$, respectively, which are valid. On $D_3^o$, the inequality is verified easily. On $D_4^o$, the estimates $U_{xx}+U_{yy}\le 0$, $U_{xx}-U_{yy}\le 0$ are equivalent to

$$2(p-2)(1-x-y)\le 0 \quad \text{and} \quad -4(1-y+x)\le 0,$$

respectively, which hold true. On $D_5^o$, it is obvious that $U_{xx} \le U_{yy}$, while the inequality $U_{xx} \le -U_{yy}$ reduces to

$$\frac{(h'(x+y)-1)(-h(x+y)+x-y+1)}{h(x+y)-(x+y)+1} \le 0,$$

which is guaranteed by (4.12) and the definition of $D_5$. Finally, assume $(x,y) \in D_6^o$. Then, by (4.7),

$$\frac{p^{p+1}}{2^p}(x-y+1)^{p-2} = \frac{2G'(x-y+1)}{(2+x-y-G(x-y+1))^2}$$

and it is easy to see that $U_{xx}(x,y) \le U_{yy}(x,y)$. The inequality $U_{xx}(x,y) \le -U_{yy}(x,y)$ is equivalent to

$$\frac{4(1-G'(x-y+1))}{(2+x-y-G(x-y+1))^3}(x+y-G(x-y+1)) \le 0.$$

To prove its validity, note that $G' \ge 1$, which is a consequence of (4.12), and $x+y \ge G(x-y+1)$, which is equivalent to $h(x+y) \ge x-y+1$, one of the inequalities in the definition of $D_6$.

The proof is complete. □

As in the proof of the case $p<1$, we extend the partial derivatives of the special function to the whole $\mathbb{R}_+ \times \mathbb{R}$. Let $\phi$, $\psi : \mathbb{R}_+ \times \mathbb{R}_+ \to \mathbb{R}$ be given by

$$\phi(x,y) = \begin{cases} U_x(x,y), & \text{if } (x,y) \notin \partial D_3 \cap \partial D_4, \\ U_x(x+,y), & \text{if } (x,y) \in \partial D_3 \cap \partial D_4, \end{cases}$$

$$\psi(x,y) = \begin{cases} U_y(x,y), & \text{if } (x,y) \notin \partial D_0 \cup (\partial D_3 \cap \partial D_4), \\ U_y(x+,y), & \text{if } (x,y) \in \partial D_3 \cap \partial D_4 \setminus \{(0,0)\}, \\ U_y(x,y+), & \text{if } (x,y) \in \partial D_0 \end{cases}$$

and extend them to the whole $\mathbb{R} \times \mathbb{R}$ by $\phi(x,y) = \phi(x,-y)$, $\psi(x,y) = -\psi(x,-y)$.

The further properties of $U$ are described in the following lemma:

**Lemma 5.4.** (i) *Let $x \ge 0, x+h \ge 0$ and $y,k \in \mathbb{R}$. Then*

$$U(x+h,y+k) \le U(x,y) + \phi(x,y)h + \psi(x,y)k. \tag{5.5}$$

(ii) *Let $x \in \mathbb{R}_+$ and $y \in (-1,1)$. Then the function $H_{x,y}$, defined on $\{t : x+t \ge 0$ and $-1 < y+t < 1\}$ and given by*

$$H_{x,y}(t) = \phi(x+t,y+t) - \psi(x+t,y+t),$$

*is non-increasing.*

**Proof.** (i) Consider a continuous function $L = L_{x,y,h,k}$ defined on $\{t: x+th \geq 0\}$ and given by

$$L(t) = u(x+th, y+tk).$$

The inequality (5.5) is equivalent to $L(1) \leq L(0) + L'(0)$ (with $L'(0)$ replaced by a left- or right-sided derivative if $(x,y)$ belongs to $\partial D_0$ or $\partial D_3 \cap \partial D_4$) and will follow if we show that $L$ is concave. To this end, it suffices to prove that $L''(t) \leq 0$ for those $t$, where the second derivative exists, and $L'(t-) \geq L'(t+)$ for remaining $t$, satisfying $x+th>0$. The first inequality follows from

$$L''(t) = U_{xx}(x+th, y+tk)h^2 + 2U_{xy}(x+th, y+tk)hk + U_{yy}(x+th, y+tk)k^2 \leq 0,$$

due to (5.1). To deal with the second, recall that $L$ is of class $C^1$ except for $(x+th, t+tk)$ belonging to $\partial D_0$ or $\partial D_3 \cap \partial D_4$. Hence, by the transity property $L_{x,y,h,k}(t+s) = L_{x+th,y+tk,h,k}(s)$, all we need is $L'(0-) \geq L'(0+)$ for $(x,y)$ belonging to one of these sets. As $L_{x,y,h,k}$ is concave if and only if $L_{x,y,-h,-k}$ is concave, we may also assume $h>0$. Now, if $(x,y) \in \partial D_0$, then $U_x$ is continuous in $(x,y)$ and

$$L'(0-) - L'(0+) = (U_y(x,y-) - U_y(x,y+))|k| = U_y(x,y-)|k| \geq 0$$

by Lemma 5.2. Suppose then that $(x,y) \in \partial D_3 \cap \partial D_4$. We have

$$L'(0-) - L'(0+) = \left(-\frac{p^2}{2(p-1)} + \frac{1}{2x} + \frac{p^p}{2(p-1)}x^{p-1}\right)(h+k) \geq 0. \tag{5.6}$$

The latter inequality is a consequence of $h \geq k$ and

$$-\frac{p^2}{2(p-1)} + \frac{1}{2x} + \frac{p^p}{2(p-1)}x^{p-1} = \frac{1}{2x(p-1)}[p(1-px) - (1-(px)^p)] \geq 0,$$

which follows from the mean value property.

(ii) If $(x+t, y+t)$ lies in the interior of one of the sets $D_k$, $k = 1,2,\ldots,7$, then

$$H'_{x,y}(t) = U_{xx}(x+t, y+t) - U_{yy}(x+t, y+t) \leq 0.$$

Therefore we will be done if we show that $H_{x,y}$ is continuous. By Lemma 5.1, we only need to check continuity for $t$ determined by $(x+t, y+t) \in \partial D_3 \cap \partial D_4$, $(x+t, y+t) \neq (2/p, \pm(p-2)/p)$. If $y+t>0$, then one can check, using the formulae for $U_x$, $U_y$, that

$$\lim_{t' \to t} H_{x,y}(t') = \lim_{t' \to t}[U_x(x+t', y+t') - U_y(x+t', y+t')] = -\frac{p^2}{2(p-1)} = H_{x,y}(t).$$

For $y+t<0$ this follows from the fact that $H_{x,y}(s) = U_x(x+s+, y+s) - U_y(x+s+, y+s)$ for $s$ lying in some neighbourhood of $t$ (both the partial derivatives are defined by the formulae for $D_4$ and hence are continuous). □

**Lemma 5.5.** (i) *For any $x \geq 0$, $y \in \mathbb{R}$ satisfying $|y| \leq x$ we have $U(x,y) \leq 0$.*
(ii) *We have $U \geq V$.*

**Proof.** (i) Using (5.5), we may write

$$U(x,y) \leq U(0,0) + \phi(0,0)x + \psi(0,0)y = 0,$$

as claimed.

(ii) The inequality is clear on $D_0$. For $(x,y) \notin D_0$, use Lemma 5.2 to obtain

$$U(x,y) \geq U(x,0) = V(x,0) = V(x,y).$$

This finishes the proof. □

## 6. The proofs of the inequalities (1.6) and (1.8)

For the sake of convenience, the proof is divided into a few steps.

*Step* 1. We start with a smoothing argument and correct the function $U$ in such a way that the key properties are still valid (the inequalities (6.1), (6.2) and (6.3) below). Let $\varepsilon > 0$ be fixed and $m$ be a positive integer satisfying $1/m < \varepsilon$. Let $g_m : \mathbb{R}^2 \to \mathbb{R}_+$ be a $C^\infty$ function with support inside the ball $B_m$ centered at 0 and radius $1/m$. Furthermore, assume $g_m$ has integral 1 and define $U^m : [1/m, \infty) \times \mathbb{R} \to \mathbb{R}$ as the convolution of $U$ with the function $g_m$. Note that $U^m$ is infinitely differentiable. Moreover, by Lemma 5.4 (i), $U^m$ is concave along the lines of a slope not greater than 1 in absolute value, as convolving with a positive function does not affect this property. Therefore, we have

$$U^m_{xx} \pm 2U^m_{xy} + U^m_{yy} \leq 0 \tag{6.1}$$

and, for $x, x+h > 1/m$, $y, k \in \mathbb{R}$ such that $|h| \geq |k|$,

$$U^m(x+h, y+k) \leq U^m(x,y) + U^m_x(x,y)h + U^m_y(x,y)k. \tag{6.2}$$

Furthermore, by part (ii) of this lemma, for any $x > 1/m$, $|y| < 1 - 1/m$, the function $H^m_{x,y}$, defined on a (small) neighbourhood of 0 by

$$H^m(t) = U^m_x(x+t, y+t) - U^m_y(x+t, y+t),$$

is non-increasing. To see this, note that by integration by parts and continuity of $U$ on $(0,\infty) \times \mathbb{R}$, we have

$$H^m_{x,y}(t) = \int_{B_m} H_{x-u,y-v}(t) g_m(u,v)\, \mathrm{d}u\, \mathrm{d}v$$

(here, $H$ under the integral is defined as in Lemma 5.4 (ii)). Thus

$$(H^m_{x,y})'(0) = U^m_{xx}(x,y) - U^m_{yy}(x,y) \leq 0 \qquad \text{for } x > 1/m,\ |y| < 1 - 1/m. \tag{6.3}$$

*Step* 2. Here we will introduce the key stopping time. With no loss of generality we may assume $\|M\|_p < \infty$. For $\eta < 1$, let $R = R(\eta, \varepsilon)$ denote the greatest number $r$ such that

$$U(x,y) \geq \eta 1_{\{|y|\geq 1-\varepsilon\}} - \frac{p^p}{2^p(p-1)} x^p \qquad \text{for all } x \leq r,\ y \in \mathbb{R}. \tag{6.4}$$

Note that if $\eta$ is fixed and $\varepsilon \downarrow 0$, then $R(\eta, \varepsilon) \to \infty$. This is a consequence of Lemma 5.5(ii).

Let $\tilde{M}_t = M_t + \varepsilon$ and introduce the stopping time

$$\tau = \inf\{t : |N_t| \geq 1 - \varepsilon \text{ or } \tilde{M} \geq R\}.$$

*Step* 3. We have

$$\mathbb{E}U^m(\tilde{M}_{\tau\wedge t}, N_{\tau\wedge t}) \leq \mathbb{E}U^m(\tilde{M}_0, N_0). \tag{6.5}$$

This can be proved essentially in the same manner as (2.5): apply Itô's formula, group the expressions under the integrals in appropriate way and use inequalities (6.1), (6.2) and (6.3) together with the differential subordination of $N$ by $M$.

*Step* 4. This is the final part. Note that we have the estimate

$$1 \geq U^m(x,y) \geq \inf_{(x',y')\in B_m} U(x-x', y-y') \geq -\frac{p^p}{2^p(p-1)}(x+\varepsilon)^p. \tag{6.6}$$

The martingales $(\tilde{M}_{\tau\wedge t}), (N_{\tau\wedge t})$ converge almost surely, so by Lebesgue's dominated convergence theorem, if we let $t \to \infty$ in (6.5), we get

$$\mathbb{E}U^m(\tilde{M}_\tau, N_\tau) \leq \mathbb{E}U^m(\tilde{M}_0, N_0).$$

We have $U^m \to U$ pointwise as $m \to \infty$ and by (6.6) we may again use Lebesgue's theorem to obtain

$$\mathbb{E}U(\tilde{M}_0, N_0) \geq \mathbb{E}U(\tilde{M}_\tau, N_\tau) = \mathbb{E}U(\tilde{M}_\tau, N_\tau)1_{\{M_\tau < R\}} + \mathbb{E}U(\tilde{M}_\tau, N_\tau)1_{\{M_\tau \geq R\}}.$$

Now $\mathbb{E}U(\tilde{M}_0, N_0) \leq 0$ by Lemma 5.5(i). For the first expression appearing on the right, we use the inequality (6.4). For the second one, we use the bound

$$U(x,y) \geq -\frac{p^p}{2^p(p-1)} x^p,$$

which is a trivial consequence of Lemma 5.5(ii). We arrive at

$$0 \geq \eta \mathbb{P}(|N_\tau| \geq 1 - \varepsilon, \tilde{M}_\tau < R) - \frac{p^p}{2^p(p-1)} \mathbb{E}\tilde{M}_\tau^p.$$

Therefore,

$$\eta \mathbb{P}(N^* \geq 1) \leq \eta \mathbb{P}(N^* \geq 1, \tilde{M}_\tau < R) + \eta \mathbb{P}(\tilde{M}_\tau \geq R)$$

$$\leq \eta\mathbb{P}(N_\tau \geq 1-\varepsilon, \tilde{M}_\tau < R) + \eta\mathbb{P}(\tilde{M}_\tau \geq R)$$
$$\leq \frac{p^p}{2^p(p-1)}\mathbb{E}\tilde{M}_\tau^p + \frac{\eta}{R^p}\mathbb{E}\tilde{M}_\tau^p \leq \left(\frac{p^p}{2^p(p-1)} + \frac{\eta}{R^p}\right)\|\tilde{M}\|_p^p$$
$$\leq \left(\frac{p^p}{2^p(p-1)} + \frac{\eta}{R^p}\right)(\varepsilon + \|M\|_p)^p.$$

Here, in the third passage, we used Chebyshev's inequality. Now let $\varepsilon \to 0$ to obtain

$$\eta\mathbb{P}(N^* \geq 1) \leq \frac{p^p}{2^p(p-1)}\|M\|_p^p.$$

As $\eta$ was arbitrary, the proof is complete.

## 7. The proof of Theorem 1.6

Let us begin with the following inequality. We omit the straightforward proof.

**Lemma 7.1.** *Let $1 \leq p \leq 2$ and $x, h \in \mathbb{R}$. Then*

$$|x+h|^p + |x-h|^p \leq 2|x|^p + 2|h|^p. \tag{7.1}$$

Now, let $(B, \mathbb{P}) = ((B^1, B^2), (\mathbb{P}_{x,y})_{(x,y)\in\mathbb{R}^2})$ denote the family of two-dimensional Brownian motions such that for any $(x, y)$,

$$\mathbb{P}_{x,y}(B_0 = (x, y)) = 1.$$

Recall the constant $K_p$ given by (1.10) and define $V : \mathbb{R}^2 \to \mathbb{R}$ by $V(x, y) = 1_{\{|y|\geq 1\}} - K_p^p|x|^p$. Throughout this section, $\tau$ will denote the stopping time

$$\tau = \inf\{t \geq 0 : |B_t^2| \geq 1\}.$$

The special function is defined by

$$U(x, y) = \mathbb{E}_{x,y}|B_\tau^1|^p.$$

Note that we have $U(x, y) = |x|^p$ for $|y| \geq 1$. It follows from the very definition that $U$ is harmonic on $S$ and continuous on $\mathbb{R}^2$. Moreover, it is obvious that for any $x$, $y$ we have

$$U(x, y) = U(-x, y) = U(x, -y) = U(-x, -y). \tag{7.2}$$

We study the further properties of the function $U$ in the lemmas below. First we will provide an explicit formula for $U$. To do this, let $H = \{(\alpha, \beta) : \beta > 0\}$ denote the upper half-plane and define $\mathcal{W} : H \to \mathbb{R}$ as the Poisson integral

$$\mathcal{W}(\alpha, \beta) = \frac{2^p}{\pi^{p+1}}\int_{-\infty}^{\infty} \frac{\beta|\log|t||^p}{(\alpha-t)^2 + \beta^2}\,dt.$$

Clearly, $\mathcal{W}$ is harmonic on $H$. Furthermore, we have

$$\lim_{(\alpha,\beta)\to(t,0)} \mathcal{W}(\alpha,\beta) = \left(\frac{2}{\pi}\right)^p |\log|t||^p. \tag{7.3}$$

Consider the conformal map $\phi$ on $S = \{(x,y): |y| < 1\}$, defined by

$$\phi(x,y) = \phi(z) = \mathrm{i}\mathrm{e}^{\pi z/2}.$$

It is easy to check that $\phi$ maps $S$ onto $H$.

**Lemma 7.2.** (i) *We have the identity*

$$U(x,y) = \begin{cases} |x|, & \text{if } |y| \geq 1, \\ \mathcal{W}(\phi(x,y)), & \text{if } |y| < 1. \end{cases} \tag{7.4}$$

(ii) *We have* $U(0,0) = K_p^{-p}$.

**Proof.** (i) Note that both sides of (7.4) are harmonic on $S$, continuous on $\mathbb{R}^2$ (for the right-hand side, use (7.3)) and equal on the set $\{(x,y): |y| \geq 1\}$. This proves the claim.

(ii) Using the first part of the lemma, we may write

$$\begin{aligned} U(0,0) = \mathcal{W}(0,1) &= \left(\frac{2}{\pi}\right)^{p+1} \int_0^\infty \frac{|\log t|^p}{t^2+1}\,\mathrm{d}t \\ &= \left(\frac{2}{\pi}\right)^{p+1} \int_{-\infty}^\infty \frac{|s|^p \mathrm{e}^s}{\mathrm{e}^{2s}+1}\,\mathrm{d}s \\ &= 2\left(\frac{2}{\pi}\right)^{p+1} \int_0^\infty s^p \mathrm{e}^{-s} \sum_{k=0}^\infty (-\mathrm{e}^{-2s})^k\,\mathrm{d}s \\ &= 2\left(\frac{2}{\pi}\right)^{p+1} \Gamma(p+1) \sum_{k=0}^\infty \frac{(-1)^k}{(2k+1)^{p+1}} = K_p^{-p}. \end{aligned}$$

Here we have used the identity

$$\sum_{k=0}^\infty \frac{1}{(2k+1)^2} = \frac{\pi^2}{8}.$$

□

**Lemma 7.3.** *For fixed* $x \in \mathbb{R}$, *the function* $U(x,\cdot)$ *is concave on* $(-1,1)$.

**Proof.** Since $U$ is harmonic on $S$, an equivalent formulation is that for any fixed $y \in (-1,1)$ the function $U(\cdot,y)$ is convex. However, the estimate

$$\lambda U(x_1,y) + (1-\lambda)U(x_2,y) \geq U(\lambda x_1 + (1-\lambda)x_2, y), \qquad \lambda \in (0,1),$$

follows directly from the definition of $U$ and Jensen's inequality. □

**Lemma 7.4.** (i) *For any* $x>0, y\in(0,1)$, $U_{xyy}\geq 0$.
(ii) $U_{xy}\geq 0$ *on* $S$.

**Proof.** (i) Since $U_x$ is harmonic on $S$, the equivalent statement is $U_{xxx}\leq 0$. Fix $x,y$ as in the statement and let $\varepsilon\in(0,x)$. We have

$$2U_x(x,y)-U_x(x-\varepsilon,y)-U_x(x+\varepsilon,y)=\mathbb{E}_{0,y}f(B^1_\tau),$$

where

$$f(h)=2|x+h|^{p-2}(x+h)-|x-\varepsilon+h|^{p-2}(x-\varepsilon+h)-|x+\varepsilon+h|^{p-2}(x+\varepsilon+h).$$

Note that $f\geq 0$ for $h\geq -x$ and $f(-x+h)=-f(-x-h)$ for any $h$. As $\tau$ is independent of $B^1$, we infer that the density of $B^1_\tau$ under $\mathbb{P}_{0,y}$ is decreasing on $[0,\infty)$. This implies $\mathbb{E}_{0,y}f(B^1_\tau)\geq 0$ and, since $x>0$ and $\varepsilon\in(0,x)$ were arbitrary, we conclude that $U_x$ is concave.

(ii) We have $U_y(x,0)=0$ due to (7.2). Thus $U_{xy}(x,0)=0$ and the claim follows from the first part of the lemma. □

**Lemma 7.5.** (i) *We have* $U(x,y)\geq U(0,0)$ *for any* $x$, $y$ *such that* $|y|\leq|x|$.
(ii) *We have* $U(x,y)\leq |x|^p+K_p^{-p}1_{\{|y|<1\}}$ *for any* $x,y\in\mathbb{R}$.

**Proof.** (i) For fixed $x$, the function $U(x,\cdot)$ is even, concave on $(-1,1)$ and constant on $(-\infty,1]\cup[1,\infty)$. This implies that it suffices to show the inequality for $x=y$. However, it is easy to check that the function $F$, given by $F(x)=U(x,x)$, $x\in\mathbb{R}$, is even, convex on $[-1,1]$ and increasing on $[1,\infty)$. This yields the claim.

(ii) It suffices to prove the inequality on $S$. The function

$$y\mapsto U(x,y)-|x|^p-K_p^{-p}1_{\{|y|<1\}}=U(x,y)-|x|^p-K_p^{-p}$$

is even (by (7.2)) and concave (due to Lemma 7.3), so attains its maximum at 0. Hence we must show that

$$U(x,0)-|x|^p-U(0,0)\leq 0,$$

or

$$\mathbb{E}_{0,0}|x+B^1_\tau|^p\leq |x|^p+\mathbb{E}_{0,0}|B^1_\tau|^p,$$

an inequality which follows from (7.1). □

**Proof of (1.6).** We may assume $\|M\|_p<\infty$. If $N$ is differentially subordinate and orthogonal to $M$, then $N$ has continuous paths (see Lemma 2.1 in Bañuelos and Wang (2000)). We proceed as in the proof of Theorem 1.5. Introduce the stopping time

$$\sigma=\inf\{t:|N_t|\geq 1\}.$$

Applying Itô's formula and using Lemma 7.3 together with the differential subordination and the orthogonality property, we get (see the proof of (2.5))

$$\mathbb{E}U(M_{\sigma\wedge t}, N_{\sigma\wedge t}) \geq \mathbb{E}U(M_0, N_0).$$

This, combined with Lemma 7.2 (ii) and Lemma 7.5, leads to

$$\mathbb{E}|M_{\sigma\wedge t}|^p + K_p^{-p}\mathbb{P}(|N_{\sigma\wedge t}| < 1) \geq K_p^{-p},$$

or

$$\mathbb{P}(|N_{\sigma\wedge t}| \geq 1) \leq K_p^p \|M\|_p^p.$$

Now the proof is completed using the scaled martingales $M^\varepsilon$, $N^\varepsilon$ in exactly the same manner as in the proof the inequality (1.7). □

**Sharpness.** Let $M = (B^1_{\tau\wedge t})$, $N = (B^2_{\tau\wedge t})$. Then $N$ is differentially subordinate and orthogonal to $M$. Moreover,

$$\mathbb{P}(N^* \geq 1) - K_p^p\|M\|_p = 1 - K_p^p U(0,0) = 0,$$

which implies that the best constant in (1.6) is not smaller than $K_p$. □

## 8. Inequalities for harmonic functions

In this section we study related weak-type inequalities for harmonic functions on Euclidean domains. Let $n$ be a positive integer and let $D$ be an open connected subset of $\mathbb{R}^n$. Fix a point $\xi$, which belongs to $D$. For two real-valued harmonic functions $u$, $v$ on $D$, we say that $v$ is *differentially subordinate* to $u$ if the following two conditions are satisfied:

$$|v(\xi)| \leq |u(\xi)| \tag{8.1}$$

and

$$|\nabla v(x)| \leq |\nabla u(x)| \qquad \text{for any } x \in D. \tag{8.2}$$

This concept was introduced by Burkholder (1989); see this paper for more information and references. We say that $u$, $v$ are orthogonal if

$$\nabla v(x) \cdot \nabla u(x) = 0 \qquad \text{for any } x \in D, \tag{8.3}$$

where $\cdot$ denotes the scalar product in $\mathbb{R}^n$.

Let $D_0$ be a bounded domain satisfying $\xi \in D_0 \subset D_0 \cup \partial D_0 \subset D$. Let $\mu^\xi_{D_0}$ denote the harmonic measure on $\partial D_0$, corresponding to $\xi$. Consider weak and strong $p$-norms given

by

$$\|u\|_{p,\infty} = \sup\sup_{\lambda}[\lambda^p \mu^{\xi}_{D_0}(\{x \in \partial D_0 : |u(x)| \geq \lambda\})]^{1/p},$$

$$\|u\|_p = \sup\left[\int_{\partial D_0} |u(x)|^p \, \mathrm{d}\mu^{\xi}_{D_0}(x)\right]^{1/p},$$

the first supremums being taken over all subdomains $D_0$ as above.

Now we state the harmonic analogue of Theorems 1.4 and 1.5.

**Theorem 8.1.** *Suppose $v$ is differentially subordinate to $u$ and $u$ is non-negative.*

(i) *For $0 < p < 1$ we have*

$$\|v\|_{p,\infty} \leq 2\|u\|_p$$

*and the inequality is sharp.*

(ii) *For $p \geq 2$ we have*

$$\|v\|_{p,\infty} \leq \frac{p}{2}(p-1)^{-1/p}\|u\|_p.$$

We do not know if the constant in (ii) is the best possible (except for the case $p = 2$, where, clearly, it is).

**Proof.** It suffices to show that for any bounded domain $D_0$ as above we have

$$\mu^{\xi}_{D_0}(\{x \in \partial D_0 : |v(x)| \geq 1\}) \leq \frac{p^p}{2^p(p-1)} \int_{\partial D_0} u(x)^p \, \mathrm{d}\mu^{\xi}_{D_0}(x). \tag{8.4}$$

Let $B = (B_t)_{t \geq 0}$ be a Brownian motion in $\mathbb{R}^n$ starting from $\xi$ and introduce a stopping time

$$\tau = \tau_{D_0} = \inf\{t : B_t \in \partial D_0\}.$$

Let $M$, $N$ be martingales defined by

$$M_t = u(B_{\tau\wedge t}) \quad \text{and} \quad N_t = v(B_{\tau\wedge t}), \qquad t \geq 0. \tag{8.5}$$

Note that $M$ is non-negative and $N$ is real-valued. The property (8.2) gives $|N_0| \leq |M_0|$ and, combined with (8.1), implies that $N$ is differentially subordinate to $M$: this follows from identities

$$[M, M]_t = |M_0|^2 + \int_0^{\tau\wedge t} |\nabla u(B_s)|^2 \, \mathrm{d}s,$$

$$[N, N]_t = |N_0|^2 + \int_0^{\tau\wedge t} |\nabla v(B_s)|^2 \, \mathrm{d}s.$$

Hence

$$\mu_{D_0}^{\xi}(\{x \in \partial D_0 : |v(x)| \geq 1\}) \leq \mathbb{P}(N^* \geq 1) \leq \frac{p}{2(p-1)^{1/p}} \|M_\tau\|_p = \frac{p}{2(p-1)^{1/p}} \|u\|_p.$$

Now we show the constant 2 is optimal in (i). Here we use the example of Burkholder (1994). Let $n=1$, $D=(-1,3)$, $\xi=0$, $u(x)=1+x$ and $v(x)=1-x$. Then $u$, $v$ are harmonic, $u$ is non-negative and $v$ is differentially subordinate to $u$. We have $\|u\|_p = u(\xi)=1$. Furthermore, for $0<\lambda<2$ we have $|v(x)|<\lambda$ if and only if $x \in (1-\lambda, 1+\lambda)$, which implies

$$\lim_{\lambda \uparrow 2} \sup_{D_0} \lambda (\mu(|v| \geq \lambda))^{1/p} = \lim_{\lambda \uparrow 2} \lambda = 2\|u\|_p.$$

Therefore we cannot replace 2 in (i) by a smaller number. □

The version of Theorem 1.6 for harmonic functions can be stated as follows:

**Theorem 8.2.** *Suppose $v$ is differentially subordinate to $u$ and $u, v$ are orthogonal. Then for $1 \leq p \leq 2$ we have*

$$\|v\|_{p,\infty} \leq K_p \|u\|_p.$$

*The inequality is sharp.*

**Proof.** The inequality is proved by the same argumentation as above, using the martingales $M$ and $N$ given by (8.5). Their orthogonality is guaranteed by

$$[M,N]_t = M_0 N_0 + \int_0^{\tau \wedge t} \nabla u(B_s) \cdot \nabla v(B_s) \, \mathrm{d}s.$$

We omit the details.

To see that the inequality is sharp, let $n=2$, $\varepsilon>0$, $D=\mathbb{R} \times (-1-\varepsilon, 1+\varepsilon)$, $\xi=(0,0)$, $u(x,y)=x$, $v(x,y)=y$. Clearly, $u$, $v$ are harmonic and orthogonal and $v$ is differentially subordinate to $u$. For $D_0=(-R,R)\times(-1,1)$ and $\mu=\mu_{D_0}^{\xi}$, we have

$$\lim_{R\to\infty} \left[ \int_{\partial D_0} |u(x)|^p \, \mathrm{d}\mu(x) \right]^{1/p} = \|B_\tau^1\|_p = \frac{1}{K_p},$$

which gives

$$\lim_{\varepsilon \downarrow 0} \|u\|_p = \frac{1}{K_p}.$$

To complete the proof, it suffices to note that

$$\lim_{R\to\infty} \mu(|v| \geq 1) = 1 \quad \text{and hence} \quad \|v\|_{p,\infty} \geq 1,$$

which implies that $K_p$ is the best possible. □

## Acknowledgement

I would like to thank the anonymous referee for comments which made the paper more concise and easy to read. Partially supported by MEiN Grant 1 PO3A 012 29 and the Foundation for Polish Science.

## References

Abramowitz, M. and Stegun, I.A., eds. (1992). *Handbook of Mathematical Functions with Formulas, Graphs and Mathematical Tables*. New York: Dover Publications. Reprint of the 1972 ed. MR1225604

Bañuelos, R. and Wang, G. (1995). Sharp inequalities for martingales with applications to the Beurling–Ahlfors and Riesz transformations. *Duke Math. J.* **80** 575–600. MR1370109

Bañuelos, R. and Wang, G. (2000). Davis's inequality for orthogonal martingales under differential subordination. *Michigan Math. J.* **47** 109–124. MR1755259

Burkholder, D.L. (1984). Boundary value problems and sharp inequalities for martingale transforms. *Ann. Probab.* **12** 647–702. MR0744226

Burkholder, D.L. (1989). Differential subordination of harmonic functions and martingales. In *Harmonic Analysis and Partial Differential Equations (El Escorial, 1987). Lecture Notes in Mathematics* **1384** 1–23. Berlin: Springer. MR1013814

Burkholder, D.L. (1991). Explorations in martingale theory and its applications. In *École d'Ete de Probabilités de Saint-Flour XIX. Lecture Notes in Mathematics* **1464** 1–66. Berlin: Springer. MR1108183

Burkholder, D.L. (1994). Strong differential subordination and stochastic integration. *Ann. Probab.* **22** 995–1025. MR1288140

Burkholder, D.L. (1999). Some extremal problems in martingale theory and harmonic analysis. In *Harmonic Analysis and Partial Differential Equations (Chicago, IL)* 99–115. Chicago, IL: Univ. Chicago Press. MR1743857

Choi, C. (1998). A weak-type inequality for differentially subordinate harmonic functions. *Trans. Amer. Math. Soc.* **350** 2687–2696. MR1617340

Dellacherie, C. and Meyer, P.A. (1982). *Probabilities and Potential B*. Amsterdam: North-Holland. MR0745449

Suh, Y. (2005). A sharp weak type $(p,p)$ inequality $(p>2)$ for martingale transforms and other subordinate martingales. *Trans. Amer. Math. Soc.* **357** 1545–1564 (electronic). MR2115376

Wang, G. (1995). Differential subordination and strong differential subordination for continuous time martingales and related sharp inequalities. *Ann. Probab.* **23** 522–551. MR1334160